\documentclass[11pt]{article}

\usepackage{latexsym,amssymb,color}
\usepackage[abbrev]{amsrefs}
\usepackage{amsmath}

\usepackage[truedimen,margin=30truemm, bottom=30truemm]{geometry}

\newtheorem{Theorem}{\bf Theorem}[section]
\newtheorem{Lemma}{\bf Lemma}[section]
\newtheorem{Proposition}{\bf Proposition}[section]
\newtheorem{Corollary}{\bf Corollary}[section]
\newtheorem{Remark}{\bf Remark}[section]
\newtheorem{Example}{\bf Example}[section]
\newtheorem{Definition}{\bf Definition}[section]

\newenvironment{theorem}{\begin{Theorem}$\!\!\!$}{\end{Theorem}}
\newenvironment{lemma}{\begin{Lemma}$\!\!\!$}{\end{Lemma}}

\newenvironment{definition}{\begin{Definition}$\!\!\!$}{\end{Definition}}

\numberwithin{equation}{section}

\begin{document}
\title{Solvability of the heat equation on a half-space\\
with a dynamical boundary condition\\ and unbounded initial data}

\author{
Marek Fila\\
Department of Applied Mathematics and Statistics,\\
Comenius University,\\
84248 Bratislava, Slovakia\\
fila@fmph.uniba.sk
\\
\\
Kazuhiro Ishige\\
Graduate School of Mathematical Sciences,\\ 
The University of Tokyo\\
3-8-1 Komaba, Meguro-ku, Tokyo 153-8914, Japan\\
ishige@ms.u-tokyo.ac.jp
\\
\\
Tatsuki Kawakami\footnote{Corresponding author, ORCID: 0000-0002-0017-2502}\\
Faculty of Advanced Science and Technology,\\ 
Ryukoku University,\\
1-5 Yokotani, Seta Oe-cho, Otsu, Shiga 520-2194, Japan.\\
kawakami@math.ryukoku.ac.jp
}
\date{}
\maketitle

\begin{abstract}
\noindent
We study the linear heat equation on a halfspace with a linear dynamical boundary condition.
We are interested in
an appropriate choice of the function space of initial functions such that the problem possesses 
a solution. It was known before that bounded initial data guarantee solvability. Here we extend that
result by showing that data from
a weighted Lebesgue space will also do so.


\end{abstract}
%
\section{Introduction}
Let $N\ge 2$ and ${\mathbb R}^N_+:={\mathbb R}^{N-1}\times{\mathbb R}_+$. 
This paper is concerned with global solvability of the problem
\begin{equation}
\label{eq:1.1}
	\left\{
	\begin{array}{ll}
	\displaystyle{\partial_t u-\Delta u=0}, & x\in{\mathbb R}^N_+,\,\,\,t>0,\vspace{5pt}\\
	\displaystyle{\partial_tu+\partial_\nu u=0}, & x\in\partial{\mathbb R}^N_+,\,\,\, t>0,\vspace{5pt}\\
	\displaystyle{u(x,0)=\varphi(x)},\qquad & x\in{\mathbb R}^N_+,
	\vspace{5pt}\\
	\displaystyle{u(x,0)=0},\qquad & x=(x',0)\in\partial{\mathbb R}^N_+,\quad x':=(x_1, x_2, \dots, x_{N-1}),
	\end{array}
\right.
\end{equation}
where $\partial_t:=\partial/\partial t$, and $\partial_\nu:=-\partial/\partial x_N$. 
The boundary condition from
\eqref{eq:1.1} describes diffusion through the boundary in processes such as thermal contact with a perfect conductor or
diffusion of solute from a well-stirred fluid or vapour (see e.g. \cite{C}). Various aspects of analysis of parabolic
equations with dynamical boundary conditions have been treated by many authors, see for example
\cites{AQR, BBR, BC, BP1, BP2, BPR, EMR, E2, E3, FIK08, FIKL2, FQ2, FV, GM, H, KI, VV1, VV2, V1}. 

In this paper we focus on the simplest linear
problem from a point of view which has not been considered yet (as far as we know). Namely, we are interested in
an appropriate choice of the function space of initial functions $\varphi$ such that problem~\eqref{eq:1.1} is solvable.

Throughout this paper we often identify ${\mathbb R}^{N-1}$ with $\partial{\mathbb R}^N_+$. 
We introduce some notation. 
Let $\Gamma_D=\Gamma_D(x,y,t)$ be the Dirichlet heat kernel on ${\mathbb R}^N_+$, that is,
\begin{equation}
\label{eq:1.2}
	\Gamma_D(x,y,t):=(4\pi t)^{-\frac{N}{2}}
	\left[\exp\biggr(-\frac{|x-y|^2}{4t}\biggr)-\exp\biggr(-\frac{|x-y_*|^2}{4t}\biggr)\right]
\end{equation}
for $(x,y,t)\in\overline{{\mathbb R}^{N}_+}\times{\mathbb R}^{N}_+\times(0,\infty)$,
where $y_*=(y',-y_N)$ for $y=(y',y_N)\in{\mathbb R}^N_+$.
Define
\begin{equation}
\label{eq:1.3}
	[S_1(t)\phi](x):=\int_{{\mathbb R}^N_+}\Gamma_D(x,y,t)\phi(y)\, dy,
	\quad (x,t)\in{\mathbb R}^N_+\times(0,\infty),
\end{equation}
for any measurable function $\phi$ in ${\mathbb R}^N_+$. 
For $x=(x',x_N)\in\overline{{\mathbb R}^N_+}$ and $t>0$, 
set
\begin{equation}
\label{eq:1.4}
	P(x',x_N,t):=C_N(x_N+t)^{1-N}\left(1+\left|\frac{x'}{x_N+t}\right|^2\right)^{-\frac{N}{2}},
\end{equation}
where $C_N$ is the constant chosen so that 
$$
	\int_{{\mathbb R}^{N-1}}P(x',x_N,t)dx'=1\quad\mbox{for all $x_N\ge 0$ and $t>0$}. 
$$
Then $P=P(x',x_N,t)$ is the fundamental solution of 
the Laplace equation in ${\mathbb R}^N_+$ 
with the homogeneous dynamical boundary condition, that is, 
$P$ satisfies 
$$
	\left\{
	\begin{array}{ll}
	\displaystyle{-\Delta P=0}, & x\in{\mathbb R}^N_+,\,\,\,t>0,\vspace{5pt}\\
	\displaystyle{\partial_tP+\partial_\nu P=0}, & x\in\partial{\mathbb R}^N_+,\,\,\, t>0,\vspace{5pt}\\
	\displaystyle{P(x,0)=\delta(x')},\qquad & x=(x',0)\in\partial{\mathbb R}^N_+,
	\end{array}
	\right.
$$
where $\delta=\delta(\cdot)$ is the Dirac delta function on $\partial{\mathbb R}^N_+={\mathbb R}^{N-1}$.  
Define
\begin{equation}
\label{eq:1.5}
	[S_2(t)\psi](x):=\int_{{\mathbb R}^{N-1}}P(x'-y',x_N,t)\psi(y')\, dy',
	\quad (x,t)\in{\mathbb R}^N_+\times(0,\infty),
\end{equation}
for any measurable function $\psi$ in ${\mathbb R}^{N-1}$. 
\vspace{5pt}

Consider
\begin{equation}
\label{eq:1.6}
	\left\{
	\begin{array}{ll}
	\displaystyle{\partial_t v=\Delta v-F[v]},
	\qquad & x\in{\mathbb R}^N_+,\,\,\,t>0,\vspace{5pt}\\
	\displaystyle{\Delta w=0}, & x\in{\mathbb R}^N_+,\,\,\,t>0,\vspace{5pt}\\
	\displaystyle{v=0},\quad
	\displaystyle{\partial_tw-\partial_{x_N} w=\partial_{x_N}v}, & x\in\partial{\mathbb R}^N_+,\,\,\, t>0,\vspace{5pt}\\
	\displaystyle{v(x,0)=\varphi(x)},\qquad & x\in{\mathbb R}^N_+,\vspace{5pt}\\
	\displaystyle{w(x,0)=0}, & x=(x',0)\in\partial{\mathbb R}^N_+,
	\end{array}
	\right.
\end{equation}
where
\begin{equation}
\label{eq:1.7}
	\begin{aligned}
	F[v](x,t) 
	&
	:=\int_{{\mathbb R}^{N-1}}P(x'-y',x_N,0)\partial_{x_N}v(y',0,t)\,dy'\\
 	&
	\qquad
	+\int_0^t\int_{{\mathbb R}^{N-1}}\partial_t P(x'-y',x_N,t-s)\partial_{x_N}v(y',0,s)\,dy'\,ds. 
	\end{aligned}
\end{equation}
Following \cite{FIK08},
we formulate the definition of a solution of \eqref{eq:1.1}. 
\begin{definition}
\label{Definition:1.1}
	Let $\varphi$ be measurable function in ${\mathbb R}^N_+$.
	Let $0<T\le\infty$ and 
	$$
		v,\,\,\partial_{x_N}v,\,\, w\in C(\overline{{\mathbb R}^N_+}\times(0,T)). 
	$$
	We call $(v,w)$ a solution of \eqref{eq:1.6} in ${\mathbb R}^N_+\times(0,T)$
	if $v$ and $w$ satisfy
	\begin{equation*}
		\begin{aligned}
		&
		v(x,t)
		=[S_1(t)\varphi](x)-\int_0^t[S_1(t-s)F[v](s)](x)\,ds,
		\\
		&
		w(x,t)
		=\int_0^t[S_2(t-s)\partial_{x_N}v(s)](x)\, ds,
	\end{aligned}
\end{equation*}
for $x\in\overline{{\mathbb R}^N_+}$ and $t\in(0,T)$.
Then we say that $u:=v+w$ is a solution of \eqref{eq:1.1} in ${\mathbb R}^N_+\times(0,T)$. 
In the case of $T=\infty$, we call $(v,w)$ a global-in-time solution of \eqref{eq:1.6}
and $u$ a global-in-time solution of \eqref{eq:1.1}.
\end{definition} 
We are ready to state the main results of this paper.
For $1\le r\le\infty$, we write 
$|\cdot|_{L^r}:=\|\cdot\|_{L^r(\partial{\mathbb R}^N_+)}$ 
and $\|\cdot\|_{L^r}:=\|\cdot\|_{L^r({{\mathbb R}^N_+})}$ 
for simplicity. 
Furthermore, for $1\le r\le\infty$ and $\alpha\ge0$, 
we define
$$
	L^r_\alpha:=\{f\in L^r(\mathbb R^N_+) : \|f\|_{L^r_\alpha}<\infty\},
$$
where
\begin{equation}
\label{eq:1.8}
	\|f\|_{L^r_\alpha}:=
	\left\{
	\begin{array}{ll}
	\displaystyle{\bigg(\int_{\mathbb R^N_+}|f(x)|^r h(x_N)^{-\alpha r}\,dx\bigg)^{\frac{1}{r}}}
	&
	\mbox{if}\quad 1\le r<\infty,
	\vspace{5pt}\\
	\displaystyle{\|f\|_{L^\infty}}
	&
	\mbox{if}\quad r=\infty,
	\end{array}
	\right.
	\,\,\,\mbox{with}\,\,\, h(x_N):=\frac{x_N}{x_N+1}.
\end{equation}
Then we can easily show that $\|f\|_{L^r_\alpha}\le \|f\|_{L^r_\beta}$ for $r\in[1,\infty]$ and $0\le\alpha\le\beta$.

\begin{theorem}
\label{Theorem:1.1}
	Let $N\ge 2$ and $1\le q\le\infty$.
	Furthermore, let
	$$
		p\in(Nq/(N-1),\infty]\quad\mbox{if}\quad q<\infty\qquad\mbox{and}\qquad
		p=\infty\quad\mbox{if}\quad q=\infty.
	$$
	For $r\in[q,\infty]$, put
	\begin{equation}
	\label{eq:1.9}
		\alpha(r)=(N-1)\left(\frac{1}{q}-\frac{1}{r}\right)+\frac{1}{q}.
	\end{equation}
	Assume $\varphi\in L^q_{\alpha(p)}$. 
	Then problem~\eqref{eq:1.6} possesses a unique global-in-time solution $(v,w)$ 
	with the following property:
	For any $T>0$ there exists $C_T>0$ such that
	\begin{equation}
	\label{eq:1.10}
		\begin{aligned}
		&
		\sup_{0<t<T}\,\bigg[t^{\frac{N}{2}(\frac{1}{q}-\frac{1}{p})}
		\bigg(\|v(t)\|_{L^p}+t^{\frac{1}{2}}\|\partial_{x_N}v(t)\|_{L^p}\bigg)
		+t^{\frac{1}{2}}|\partial_{x_N}v(t)|_{L^r}\bigg]\le C_T\|\varphi\|_{L^q_{\alpha(p)}},
		\\
		&
		\sup_{0<t<T}\,\bigg[\|w(t)\|_{L^p}+|w(t)|_{L^r}\bigg]
		\le C_T\|\varphi\|_{L^q_{\alpha(p)}},
		\end{aligned}
	\end{equation}
	for $r\in[q,p]$.
	Furthermore, 
	$v$ and $w$ are bounded and smooth in $\overline{{\mathbb R}^N_+}\times I$ 
	for any bounded interval $I\subset(0,\infty)$.
\end{theorem}
Let us now explain the role of the space $L^q_{\alpha(p)}$ in our study. Let $1\le q\le\infty$
and take arbitrary functions $\Phi\in L^q(\mathbb R^{N-1})$, $\vartheta\in L^q(1,\infty)$.
Now set $\varphi(x):=\Phi(x')\Psi(x_N)$ for $x=(x',x_N)\in {\mathbb R}^N_+$, where
$$
	\Psi(x_N):=
	\left\{
	\begin{array}{ll}
	x^\lambda_N
	&
	\mbox{if}\quad 0< x_N\le 1,\quad \lambda\in {\mathbb R},
	\vspace{5pt}\\
	\vartheta(x_N)
	&
	\mbox{if}\quad x_N>1.
	\end{array}
	\right.
$$
Choose $p$ as in Theorem~\ref{Theorem:1.1}. Then it is easy to check that $\varphi\in L^q_{\alpha(p)}$
if and only if
$$
	\lambda >(N-1)\left(\frac{1}{q}-\frac{1}{p}\right)(> 0\quad \mbox{if }q<\infty).
$$
If $\lambda >0$, then $\displaystyle{\lim_{x_N\to 0}\varphi(x)}=0$
which means that the condition $u(x',0,0)=0$ in \eqref{eq:1.1} is satisfied. This indicates that the choice of the space of initial 
functions is natural and also optimal in some sense since $\lambda$ can be arbitrarily close to $0$ if $q$ is large enough.

We have not observed the importance of the behavior of $\varphi$ near $\partial{\mathbb R}^N_+$
in the $L^\infty$-setting in \cite{FIK08}. The main novelty of this paper consists in working
in an appropriate weighted $L^q$-space by which we extend a result from \cite{FIK08} significantly,
as we explain below.

In \cite{FIK08} we studied the problem
\begin{equation}
\label{eq:1.11}
	\left\{
	\begin{array}{ll}
	\displaystyle{\partial_t u-\Delta u=0}, & x\in{\mathbb R}^N_+,\,\,\,t>0,\vspace{5pt}\\
	\displaystyle{\partial_tu+\partial_\nu u=0}, & x\in\partial{\mathbb R}^N_+,\,\,\, t>0,\vspace{5pt}\\
	\displaystyle{u(x,0)=\varphi(x)},\qquad & x\in{\mathbb R}^N_+,
	\vspace{5pt}\\
	\displaystyle{u(x,0)=\varphi_b(x')},\qquad & x=(x',0)\in\partial{\mathbb R}^N_+,
	\end{array}
	\right.
\end{equation}
where $\varphi$ and $\varphi_b$ are bounded functions. 
A part of Theorem~1.1 in \cite{FIK08} reads as follows:
\begin{theorem}
\label{Theorem:1.2}
	Let $N\ge 2$, $\varphi\in L^\infty({\mathbb R}^N_+)$ and $\varphi_b\in L^\infty({\mathbb R}^{N-1})$. 
	Then problem~\eqref{eq:1.11} possesses a unique global-in-time solution 
	$u$ which is bounded and smooth in $\overline{{\mathbb R}^N_+}\times I$ 
	for any bounded interval $I\subset(0,\infty)$. 
\end{theorem}
Hence, if $\varphi_b\equiv 0$ then Theorem~\ref{Theorem:1.2} is a very special case of Theorem~\ref{Theorem:1.1}.
If $\varphi_b\in L^\infty({\mathbb R}^{N-1})$ and $\varphi\in L^q_{\alpha(p)}$ with $p, q$ as in 
Theorem~\ref{Theorem:1.1}, then we can combine Theorem~\ref{Theorem:1.1} with Theorem~\ref{Theorem:1.2}
to obtain the existence of a solution of \eqref{eq:1.11} easily, since the problem is linear.
\section{Preliminaries}

In this section
we prove several lemmata on $S_1(t)\phi$ and $F[v]$,
and recall some properties of $S_2(t)\psi$.
In what follows, by the letter $C$
we denote generic positive constants (independent of $x$ and $t$)
and they may have different values also within the same line. 
\vspace{3pt}

We first recall some properties of $S_1(t)\phi$ (see e.g., \cite{GGS} and \cite{FIK08}*{Lemma~2.1}).
\begin{itemize}
\item[($G_1$)]
	For any $1\le q\le r\le\infty$,
	$$
		\|S_1(t)\phi\|_{L^r}\le c_1 t^{-\frac{N}{2}(\frac{1}{q}-\frac{1}{r})}\|\phi\|_{L^q},\qquad t>0,
	$$
	for all $\phi\in L^q(\mathbb R^N_+)$,
	where $c_1$ is a positive constant, independent of $q$ and $r$. 
	In particular, if $q=r$, then
	$$
		\sup_{t>0}\,\|S_1(t)\phi\|_{L^r}\le \|\phi\|_{L^r}.
	$$
	Furthermore, for any $1\le q\le r\le\infty$,
	\begin{equation}
	\label{eq:2.1}
		\|\partial_{x_N}S_1(t)\phi\|_{L^r}\le c_2 t^{-\frac{N}{2}(\frac{1}{q}-\frac{1}{r})-\frac{1}{2}}\|\phi\|_{L^q},\qquad t>0,
	\end{equation}
	for all $\phi\in L^q(\mathbb R^N_+)$,
	where $c_2$ is a positive constant, independent of $q$ and $r$. 
\item[($G_2$)]
	Let $\phi\in L^q(\mathbb R^N_+)$ with $1\le q\le\infty$ and $T>0$. 
	Then $S_1(t)\phi$ is bounded and smooth with respect to $x$ and $t$ 
	in $\overline{{\mathbb R}^N_+}\times(T,\infty)$. 
\end{itemize} 
\begin{lemma}
\label{Lemma:2.1}
	Let $1\le q\le r\le \infty$. 
	Assume $\phi\in L^q_{\alpha(r)}$ with $\alpha(r)$ as in \eqref{eq:1.9}.
	Then there exists $c_3=c_3(N)>0$ such that
	\begin{equation}
	\label{eq:2.2}
	|\partial_{x_N}[S_1(t)\phi]|_{L^r}\le c_3t^{-\frac{1}{2}}\|\phi\|_{L^q_{\alpha(r)}},
	\qquad t>0.
	\end{equation}
\end{lemma}
{\bf Proof.}
Let $\Gamma_d$ $(d=1,2,\dots)$ be the Gauss kernel in $\mathbb R^d$.  
It follows from \eqref{eq:1.2} that  
\begin{align}
\label{eq:2.3}
 	& 
	\Gamma_D(x,y,t)=\Gamma_{N-1}(x'-y',t)\bigg(\Gamma_1(x_N-y_N,t)-\Gamma_1(x_N+y_N,t)\bigg),
	\\
\label{eq:2.4}
	 &
	 \begin{aligned}
 	& 
	K(x,y,t):=\partial_{x_N}\Gamma_D(x,y,t)
	\\
 	& \qquad=\Gamma_{N-1}(x'-y',t)
	\left(-\frac{x_N-y_N}{2t}\Gamma_1(x_N-y_N,t)+\frac{x_N+y_N}{2t}\Gamma_1(x_N+y_N,t)\right),
	\end{aligned}
\end{align}
for $(x,y,t)\in\overline{{\mathbb R}^N_+}\times{\mathbb R}^N_+\times(0,\infty)$.
Then we have
\begin{equation}
\label{eq:2.5}
	K(x',0,y,t) = \frac{y_N}{t}\Gamma_{N-1}(x'-y',t)\Gamma_1(y_N,t)
\end{equation}
for $(x',y,t)\in{\mathbb R}^{N-1}\times{\mathbb R}^N_+\times(0,\infty)$.

We first prove \eqref{eq:2.2} for the case $r=\infty$.
By \eqref{eq:2.5} we can easily show that, for $q\in(1,\infty]$, it holds that
\begin{equation}
\label{eq:2.6}
	\begin{aligned}
	&
	\int_{\mathbb R^N_+}
	\bigg[|K(x',0,y,t)|y_N^{\frac{N}{q}}\bigg]^{q'}\,dy
	\\
	&
	=
	t^{-\frac{q'}{2}}\int_{\mathbb R^N_+}
	\bigg[t^{\frac{N}{2q}}\Gamma_{N-1}(y',t)\bigg(\frac{y_N}{t^{1/2}}\bigg)^{1+\frac{N}{q}}
	\Gamma_1(y_N,t)\bigg]^{q'}\,dy
	\\
	&
	\le Ct^{-\frac{q'}{2}}\int_0^\infty\eta^{\left(1+\frac{N}{q}\right)q'}\exp\bigg(-C\eta^2\bigg)\,d\eta
	\le Ct^{-\frac{q'}{2}}
	\end{aligned}
\end{equation}
for $t>0$, where $1/q+1/{q'}=1$.
Furthermore, for $q=1$, it holds that
\begin{equation}
\label{eq:2.7}
	|K(x',0,y,t)|y_N^N
	=
	t^{\frac{N-1}{2}}\Gamma_{N-1}(x'-y',t)\bigg(\frac{y_N}{t^{1/2}}\bigg)^{1+N}
	\Gamma_1(y_N,t)
	\le Ct^{-\frac{1}{2}},
\end{equation}
for $(x',y,t)\in{\mathbb R}^{N-1}\times{\mathbb R}^N_+\times(0,\infty)$.
For $q\in\{1, \infty\}$, by \eqref{eq:1.3}, \eqref{eq:1.8}, \eqref{eq:2.6}, and \eqref{eq:2.7} we have
\begin{equation}
\label{eq:2.8}
	\begin{aligned}
	&
	|\partial_{x_N}[S_1(t)\phi](x',0)|
	\\
 	& 
	\le\int_{\mathbb R^N_+}|K(x',0,y,t)|
	h(y_N)^{\frac{N}{q}}h(y_N)^{-\frac{N}{q}}|\phi(y)|\,dy	
	\\
 	& 
	\le\int_{\mathbb R^N_+}|K(x',0,y,t)|
	y_N^{\frac{N}{q}}h(y_N)^{-\frac{N}{q}}|\phi(y)|\,dy	
	\le Ct^{-\frac{1}{2}}\|\phi\|_{L^q_{\alpha(\infty)}}
	\end{aligned}
\end{equation}
for $x'\in\mathbb R^{N-1}$ and $t>0$.
Furthermore, for $q\in(1,\infty)$,
similarly to \eqref{eq:2.8},
applying H\"older's inequality, we see that
$$
	\begin{aligned}
	&
	|\partial_{x_N}[S_1(t)\phi](x',0)|
	\\
 	& 
	\le C\bigg\{\int_{\mathbb R^N_+}\bigg[h(y_N)^{-\frac{N}{q}}|\phi(y)|\bigg]^q\,dy\bigg\}^{\frac{1}{q}}
	\bigg\{\int_{\mathbb R^N_+}\bigg[|K(x',0,y,t)|h(y_N)^{\frac{N}{q}}\bigg]^{q'}\,dy\bigg\}^{\frac{1}{q'}}
	\\
	&
	\le Ct^{-\frac{1}{2}}\|\phi\|_{L^q_{\alpha(\infty)}}
	\end{aligned}
$$
for $x'\in\mathbb R^{N-1}$ and $t>0$,
where $1/q+1/q'=1$.
This together with \eqref{eq:2.8} implies \eqref{eq:2.2} for the case $r=\infty$.

Next we consider the case $r<\infty$.
Let $1\le q \le r<\infty$.
Put
\begin{equation}
\label{eq:2.9}
	\frac{1}{p}:=1+\frac{1}{r}-\frac{1}{q},
	\qquad
	\beta =N\left(1-\frac{1}{p}\right).
\end{equation}
Then, for $\theta\in(0,1)$, it follows from H\"older's inequality that
\begin{equation*}
	\begin{aligned}
	&
	\int_{\mathbb R^N_+}|K(x',0,y,t)||\phi(y)|\,dy
	\\
	&
	\le \bigg\{\int_{\mathbb R^N_+}\bigg[|K(x',0,y,t)|h(y_N)^{\alpha(r)}\bigg]^p
	\bigg[h(y_N)^{-\alpha(r)}|\phi(y)|\bigg]^{(1-\theta)p}\,dy\bigg\}^{\frac{1}{p}}
	\\
	&
	\hspace{7cm}
	\times \bigg\{\int_{\mathbb R^N_+}
	\bigg[h(y_N)^{-\alpha(r)}|\phi(y)|\bigg]^{\theta p'}\,dy\bigg\}^{\frac{1}{p'}}
	\end{aligned}
\end{equation*}
for $x'\in\mathbb R^{N-1}$ and $t>0$,
where $1/p+1/p'=1$.
Put $\theta=q/p'$.
Then, since it follows from \eqref{eq:2.9} that $1-\theta=q/r$, we have
\begin{equation}
\label{eq:2.10}
	\begin{aligned}
	&
	\int_{\mathbb R^N_+}|K(x',0,y,t)||\phi(y)|\,dy
	\\
	&
	\le \bigg\{\int_{\mathbb R^N_+}\bigg[|K(x',0,y,t)|h(y_N)^{\alpha(r)}\bigg]^p
	\bigg[h(y_N)^{-\alpha(r)}|\phi(y)|\bigg]^{\frac{pq}{r}}\,dy\bigg\}^{\frac{1}{p}}
	\\
	&
	\hspace{7cm}
	\times \bigg\{\int_{\mathbb R^N_+}
	\bigg[h(y_N)^{-\alpha(r)}|\phi(y)|\bigg]^q\,dy\bigg\}^{\frac{1}{q}-\frac{1}{r}}
	\\
	&
	\le \bigg\{\int_{\mathbb R^N_+}\bigg[|K(x',0,y,t)|y_N^{\alpha(r)}\bigg]^p
	\bigg[h(y_N)^{-\alpha(r)}|\phi(y)|\bigg]^{\frac{pq}{r}}\,dy\bigg\}^{\frac{1}{p}}\|\phi\|_{L^q_{\alpha(r)}}^{1-\frac{q}{r}}
	\end{aligned}
\end{equation}
for $x'\in\mathbb R^{N-1}$ and $t>0$.
Furthermore, put $\tilde\theta=p/q'<1$.
Then, since $1-\tilde\theta=p/r$ and $\beta=N/p'$,
applying H\"older's inequality with \eqref{eq:2.6}, we see that
$$
	\begin{aligned}
	&
	\int_{\mathbb R^N_+}\bigg[|K(x',0,y,t)|y_N^{\alpha(r)}\bigg]^p
	\bigg[h(y_N)^{-\alpha(r)}|\phi(y)|\bigg]^{\frac{pq}{r}}\,dy
	\\
	&
	=\int_{\mathbb R^N_+}
	\bigg[|K(x',0,y,t)|y_N^\beta\bigg]^{\tilde\theta p}
	\bigg[|K(x',0,y,t)|
	y_N^{\frac{\alpha(r)-\tilde\theta\beta}{1-\tilde\theta}}\bigg]^{(1-\tilde\theta) p}
	\bigg[h(y_N)^{-\alpha(r)}|\phi(y)|\bigg]^{\frac{pq}{r}}\,dy
	\\
	&
	\le \bigg\{\int_{\mathbb R^N_+}
	\bigg[|K(x',0,y,t)|y_N^\beta\bigg]^p\,dy\bigg\}^{\tilde\theta}
	\\
	&
	\qquad
	\times
	\bigg\{\int_{\mathbb R^N_+}\bigg[|K(x',0,y,t)|
	y_N^{\frac{\alpha(r)-\tilde\theta\beta}{1-\tilde\theta}}\bigg]^p
	\bigg[h(y_N)^{-\alpha(r)}|\phi(y)|\bigg]^q\,dy\bigg\}^{\frac{p}{r}}
	\\
	&
	\le 
	Ct^{-\frac{p\tilde\theta}{2}}
	\bigg\{\int_{\mathbb R^N_+}\bigg[|K(x',0,y,t)|
	y_N^{\frac{\alpha(r)-\tilde\theta\beta}{1-\tilde\theta}}\bigg]^p
	\bigg[h(y_N)^{-\alpha(r)}|\phi(y)|\bigg]^q\,dy\bigg\}^{\frac{p}{r}}
	\end{aligned}
$$
for $x'\in\mathbb R^{N-1}$ and $t>0$.
This together with \eqref{eq:2.10} yields that
\begin{equation}
\label{eq:2.11}
	\begin{aligned}
	&
	\int_{\mathbb R^N_+}|K(x',0,y,t)||\phi(y)|\,dy
	\\
	&
	\le 
	Ct^{-\frac{\tilde\theta}{2}}
	\bigg\{\int_{\mathbb R^N_+}\bigg[|K(x',0,y,t)|
	y_N^{\frac{\alpha(r)-\tilde\theta\beta}{1-\tilde\theta}}\bigg]^p
	\bigg[h(y_N)^{-\alpha(r)}|\phi(y)|\bigg]^q\,dy\bigg\}^{\frac{1}{r}}\|\phi\|_{L^q_{\alpha(r)}}^{1-\frac{q}{r}}
	\end{aligned}
\end{equation}
for $x'\in\mathbb R^{N-1}$ and $t>0$.
Since
$$
	\frac{\alpha(r)-\tilde\theta\beta}{1-\tilde\theta}=(N-1)\bigg(1-\frac{1}{p}\bigg)+1,
$$
similarly to \eqref{eq:2.6}, it holds that
$$
	\begin{aligned}
	&
	\int_{\mathbb R^{N-1}}\bigg[|K(x',0,y,t)|y_N^{\frac{\alpha(r)-\tilde\theta\beta}{1-\tilde\theta}}\bigg]^p\,dx'
	\\
	&
	=t^{-\frac{p}{2}}
	\int_{\mathbb R^{N-1}}
	\bigg[t^{\frac{N-1}{2}(1-\frac{1}{p})+\frac{1}{2}}
	\Gamma_{N-1}(y',t)\bigg(\frac{y_N}{t^{1/2}}\bigg)^{(N-1)(1-\frac{1}{p})+2}
	\Gamma_1(y_N,t)\bigg]^p\,dx'
	\\
	&
	\le Ct^{-\frac{p}{2}}\bigg(\frac{y_N}{t^{1/2}}\bigg)^{(N-1)(p-1)+2p}\exp\bigg(-C\frac{y_N^2}{t}\bigg)
	\le Ct^{-\frac{p}{2}}.
	\end{aligned}
$$
Tis together with \eqref{eq:2.11} implies
$$
\begin{aligned}
	&
	|\partial_{x_N}[S_1(t)\phi]|_{L^r}^r
	\\
 	& 
	\le
	\int_{\mathbb R^{N-1}}\bigg\{\int_{\mathbb R^N_+}|K(x',0,y,t)||\phi(y)|\,dy\bigg\}^r\,dx'
	\\
	&
	\le Ct^{-\frac{r\tilde\theta}{2}}
	\int_{\mathbb R^{N-1}}\int_{\mathbb R^N_+}\bigg[|K(x',0,y,t)|
	y_N^{\frac{\alpha(r)-\tilde\theta\beta}{1-\tilde\theta}}\bigg]^p
	\bigg[h(y_N)^{-\alpha(r)}|\phi(y)|\bigg]^q\,dy\,dx'
	\|\phi\|_{L^q_{\alpha(r)}}^{r-q}
	\\
	&
	\le Ct^{-\frac{r}{2}}\|\phi\|_{L^q_{\alpha(r)}}^r
	\end{aligned}
$$
for $t>0$,
and we have \eqref{eq:2.2}.
Thus Lemma~\ref{Lemma:2.1} follows.
$\Box$ \vspace{5pt}

Next we recall some properties of $S_2(t)\psi$. 
\begin{itemize}
\item[($P_1$)] 
  	Let $\psi\in L^r({\mathbb R}^{N-1})$ for some $r\in[1,\infty]$ 
	and $t$, $t'>0$. Then
  	$$
	\begin{aligned}
    		& 
		[S_2(t)\psi](x',x_N)=[S_2(t+x_N)\psi](x',0),\vspace{3pt}
    		\\
    		& 
		[S_2(t+t')\psi](x)=[S_2(t)(S_2(t')\psi)](x),
  	\end{aligned}
	$$
	for $x=(x',x_N)\in\overline{{\mathbb R}^N_+}$. 
  	Furthermore, 
  	$$
  		\lim_{t\to 0}|S_2(t)\psi-\psi|_r=0\quad\mbox{if $1\le r<\infty$}.
  	$$
\item[($P_2$)] 
  	For any $1\le r\le q\le\infty$,  
 	$$
  		|S_2(t)\psi|_{L^q}\le Ct^{-(N-1)(\frac{1}{r}-\frac{1}{q})}|\psi|_{L^r},\qquad t>0, 
	$$
  	for all $\psi\in L^r({\mathbb R}^{N-1})$. 
  	In particular, if $q=r$, then 
	\begin{equation}
	\label{eq:2.12}
  		\sup_{t>0}\,|S_2(t)\psi|_{L^q}\le |\psi|_{L^q}.
	\end{equation}
\item[($P_3$)]
 	Let $1\le r<\infty$ and $Nr/(N-1)<q\le\infty$.
 	Then
 	$$
 	\|S_2(t)\psi\|_{L^q}\le Ct^{-(N-1)(\frac{1}{r}-\frac{1}{q})+\frac{1}{q}}|\psi|_{L^r},\qquad t>0,
 	$$
 	for all $\psi\in L^r({\mathbb R}^{N-1})$.
 	Furthermore,
 	\begin{equation}
 	\label{eq:2.13}
 		\sup_{t>0}\,\|S_2(t)\psi\|_{L^q}\le C(|\psi|_{L^q}+|\psi|_{L^r})
 	\end{equation}
 	for all $\psi\in L^q({\mathbb R}^{N-1})\cap L^r({\mathbb R}^{N-1})$. 
\end{itemize}
Properties~($P_1$), ($P_2$), and ($P_3$) easily follow from \eqref{eq:1.5} 
(see e.g. \cite{FIK02}) and imply that 
$$
	\sup_{t>0}\,\|S_2(t)\psi\|_{L^\infty}\le |\psi|_{L^\infty}
$$
for all $\psi\in L^\infty({\mathbb R}^{N-1})$. 
Furthermore,
by an argument  similar to that in the proof of property~($G_2$) 
we have: 
\begin{itemize}
\item[($P_4$)] 
  	Let $\psi\in L^r({\mathbb R}^{N-1})$ with $1\le r\le\infty$. 
  	Then, for any $T>0$, $S_2(t)\psi$ is bounded and smooth in $\overline{{\mathbb R}^N_+}\times(T,\infty)$. 
\end{itemize}

At the end of this section,
we have the following (see also \cite{FIKL2}*{Lemma~3.3}).
\begin{lemma}
\label{Lemma:2.2}
	Let $0\le a<1$ and $0\le b<1$ be such that $0\le a+b\le 1$. 
	Let $\gamma\ge 0$ and $T>0$. 
	Then, for any $\delta>0$, 
	there exists a $M_*\ge 1$ such that 
	$$
		\sup_{0<t<T}e^{-Mt}t^\gamma\int_0^t e^{Ms}s^{-a}(t-s)^{-b}\,ds\le\delta\quad\mbox{for}\quad M\ge M_*.
	$$
\end{lemma}
The proof of this lemma is almost the same as the proof of \cite{FIKL2}*{Lemma~3.3},
so we omit it here.
%

%
\section{Proof of Theorem~\ref{Theorem:1.1}}
In this section we prove Theorem~\ref{Theorem:1.1}.
By \cite{FIK08}*{Theorem~1.1} with $\varepsilon=1$ and $\varphi_b\equiv0$
we have Theorem~\ref{Theorem:1.1} for the case $p=q=\infty$.
So we focus on the case $q<\infty$.
\bigskip

Let $T>0$, $M\ge1$, $1\le q< \infty$, and $p\in(Nq/(N-1),\infty]$.
Set 
$$
	X_{T,M}:=
	\bigg\{v\,:\,v, \partial_{x_N}v\in C(\overline{{\mathbb R}^N_+}\times(0,T)), \|v\|_{X_{T,M}}<\infty\bigg\},
	\quad
	\|v\|_{X_{T,M}}:=\sup_{0<t<T}e^{-Mt}E[v](t),
$$
where 
$$
	E[v](t)
	:=
	t^{\frac{N}{2}\left(\frac{1}{q}-\frac{1}{p}\right)}\left[\|v(t)\|_{L^p}+t^{\frac{1}{2}}\|\partial_{x_N}v(t)\|_{L^p}\right]
	+\sup_{q\le r\le p}\,t^{\frac{1}{2}}|\partial_{x_N}v(t)|_{L^r}.
$$
Then $X_{T,M}$ is a Banach space equipped with the norm $\|\cdot\|_{X_{T,M}}$. 
We apply the Banach contraction mapping principle in $X_{T,M}$ 
to find a fixed point of the functional 
\begin{equation}
\label{eq:3.1}
	Q[v](t)
	:=S_1(t)\varphi-D[v](t)
\end{equation}
on $X_{T,M}$, where $D[v]$ is the function defined by
\begin{equation}
\label{eq:3.2}
	D[v](t):=\int_0^tS_1(t-s)F[v](s)\,ds
\end{equation}
and $F[v]$ is the function defined by \eqref{eq:1.7}.
\begin{lemma}
\label{Lemma:3.1}
	Let $T>0$, $M\ge1$, $1\le q<\infty$, and $p\in(Nq/(N-1),\infty]$.
	Assume that $v\in X_{T,M}$.
	Then there exists $C>0$, independent of $T$ and $M$, such that, for $p\in(Nq/(N-1),\infty)$, 
	it holds that
	\begin{equation}
	\label{eq:3.3}
		\|F[v](t)\|_{L^p}
 		\le C(1+t^{\frac{1}{p}})t^{-\frac{1}{2}}e^{Mt}\|v\|_{X_{T,M}}
	\end{equation}
	for $0<t<T$.
	Furthermore, for any $r\in[q,p]$,
	\begin{equation}
	\label{eq:3.4}
		\|F[v](\cdot,x_N,t)\|_{L^r(\mathbb R^{N-1})}
 		\le C\bigg(1+(x_N^{-1}t)^{\frac{1}{2}}\bigg)t^{-\frac{1}{2}}e^{Mt}\|v\|_{X_{T,M}}
	\end{equation}
	 for $x_N\in(0,\infty)$ and $0<t<T$.
\end{lemma}
\noindent
{\bf Proof.}
Let $T>0$, $M\ge1$, $1\le q<\infty$, $p\in(Nq/(N-1),\infty]$, and $v\in X_{T,M}$.
It follows from \eqref{eq:1.7} that 
\begin{equation}
\label{eq:3.5}
	F[v](x,t)=F_1[v](x,t)+F_2[v](x,t)
\end{equation}
for $x\in\mathbb R^N_+$ and $0<t<T$,
where 
$$
	\begin{aligned}
 	& 
	F_1[v](x,t):=\int_{{\mathbb R}^{N-1}}P(x'-y',x_N,0)\partial_{x_N}v(y',0,t)\, dy',
	\\
 	& 
	F_2[v](x,t):=\int_0^t\int_{{\mathbb R}^{N-1}}\partial_t P(x'-y',x_N,t-s)\partial_{x_N}v(y',0,s)\, dy'\,ds.
	\end{aligned}
$$
We first obtain some estimates of $F_1[v]$.
For $p\in(Nq/(N-1),\infty)$, 
by \eqref{eq:1.5} and \eqref{eq:2.13} we have 
\begin{equation}
\label{eq:3.6}
	\begin{aligned}
	\|F_1[v](t)\|_{L^p} 
	& 
	\le\liminf_{\varepsilon\to +0}
	\left(\int_{\varepsilon}^\infty\int_{{\mathbb R}^{N-1}}
	|[S_2(x_N)\partial_{x_N}v(\cdot,0,t)](x',0)|^p\,dx'\,dx_N\right)^{1/p}
	\\
 	& 
	=\liminf_{\varepsilon\to +0}\left(\int_{{\mathbb R}^N_+}|[S_2(x_N)
 	\partial_{x_N}v(\cdot,0,t)](x',\varepsilon)|^p\,dx\right)^{1/p}
	\\
	&
	=\liminf_{\varepsilon\to +0}\left(\int_{{\mathbb R}^N_+}
 	|[S_2(\epsilon)\partial_{x_N}v(\cdot,0,t)](x)|^p\,dx\right)^{1/p}
	\\
 	& 
	\le C\bigg(|\partial_{x_N}v(t)|_{L^q}+|\partial_{x_N}v(t)|_{L^p}\bigg)
	\le Ct^{-\frac{1}{2}}e^{Mt}\|v\|_{X_{T,M}}
	\end{aligned}
\end{equation}
for $0<t<T$.
Similarly, for $r\in[q,p]$, by ($P_2$) we obtain
\begin{equation}
\label{eq:3.7}
	\begin{aligned}
	\|F_1[v](\cdot,x_N,t)\|_{L^r(\mathbb R^{N-1})}
	&
	\le |[S_2(x_N)\partial_{x_N}v(\cdot,0,t)]|_{L^r}
	\\
	&
	\le |\partial_{x_N}v(t)|_{L^r}
	\le t^{-\frac{1}{2}}e^{Mt}\|v\|_{X_{T,M}}
	\end{aligned}
\end{equation}
for $x_N\in(0,\infty)$ and $0<t<T$. 

Next we obtain some estimates of $F_2[v]$. 
It follows from \eqref{eq:1.4} that
$$
	\partial_tP(x',x_N,t)
	=\frac{1}{x_N+t}\frac{|x'|^2-(N-1)(x_N+t)^2}{|x'|^2+(x_N+t)^2}P(x',x_N,t), 
$$
for $x=(x',x_N)\in\overline{{\mathbb R}^N_+}$ and $t>0$. 
This implies that
\begin{equation}
\label{eq:3.8}
	|\partial_tP(x',x_N,t)|
	\le C(x_N+t)^{-1}P(x',x_N,t)
	\le CP(x',x_N,t)x_N^{-\frac{1}{2}}t^{-\frac{1}{2}}
\end{equation}
for $x=(x',x_N)\in\overline{{\mathbb R}^N_+}$ and $t>0$. 
Furthermore, for $p\in(Nq/(N-1),\infty)$, by \eqref{eq:2.12} we have
\begin{equation}
\label{eq:3.9}
	\begin{aligned}
 	& 
	\int_{\mathbb R^N_+}\bigg\{(x_N+t-s)^{-1}[S_2(t-s)|\partial_{x_N}v(s)|](x)\bigg\}^p\,dx
	\\
 	& 
	=\int_0^\infty (x_N+t-s)^{-p}|S_2(t-s+x_N)|\partial_{x_N}v(s)|_{L^p}^p\,dx_N
	\\
 	& 
	\le C|\partial_{x_N}v(s)|_{L^p}^p\int_0^\infty (x_N+t-s)^{-p}\,dx_N
 	\le C(t-s)^{1-p}s^{-\frac{p}{2}}e^{Mps}\|v\|_{X_{T,M}}
	\end{aligned}
\end{equation}
for $0<s<t<T$. 
Then, by \eqref{eq:3.8} and \eqref{eq:3.9}  we obtain
\begin{equation}
\label{eq:3.10}
	\begin{aligned}
	&
	\|F_2[v](t)\|_{L^p}
	\\
 	& 
	\le C\int_0^t\bigg(\int_{\mathbb R^N_+}
	\bigg\{(x_N+t-s)^{-1}[S_2(t-s)|\partial_{x_N}v(s)|](x)\bigg\}^p\,dx\bigg)^{\frac{1}{p}}\,ds
	\\
 	& \le C\|v\|_{X_{T,M}}\int_0^t(t-s)^{-1+\frac{1}{p}}s^{-\frac{1}{2}}e^{Ms}\,ds
	\le Ct^{-\frac{1}{2}+\frac{1}{p}}e^{Mt}\|v\|_{X_{T,M}}
	\end{aligned}
\end{equation}
for $0<t<T$.
In addition, for $r\in[q,p]$, by ($P_2$) and \eqref{eq:3.8} we see that
\begin{equation}
\label{eq:3.11}
	\begin{aligned}
	&
	\|F_2[v](\cdot,x_N,t)\|_{L^r(\mathbb R^{N-1})} 
	\\
	&
	\le Cx_N^{-\frac{1}{2}}\int_0^t(t-s)^{-\frac{1}{2}}|[S_2(x_N+t-s)\partial_{x_N}v(\cdot,0,s)]|_{L^r}\,ds
	\\
	&
	\le Cx_N^{-\frac{1}{2}}\int_0^t(t-s)^{-\frac{1}{2}}|\partial_{x_N}v(s)|_{L^r}\,ds
	\\
	&
	\le Cx_N^{-\frac{1}{2}}\|v\|_{X_{T,M}}\int_0^t(t-s)^{-\frac{1}{2}}s^{-\frac{1}{2}}e^{Ms}\,ds
	\le Cx_N^{-\frac{1}{2}}e^{Mt}\|v\|_{X_{T,M}}
	\end{aligned}
\end{equation}
for $x_N\in(0,\infty)$ and $0<t<T.$
Therefore, by \eqref{eq:3.5}, \eqref{eq:3.6}, \eqref{eq:3.7}, \eqref{eq:3.10} and \eqref{eq:3.11} we obtain \eqref{eq:3.3} and \eqref{eq:3.4}.
Thus Lemma~\ref{Lemma:3.1} follows. 
$\Box$
\begin{lemma}
\label{Lemma:3.2}
	Assume the same conditions as in Lemma~$\ref{Lemma:3.1}$.
	Let $D[v]$ be the function defined by \eqref{eq:3.2}.
	Then there exists $M_*\ge1$ such that
	\begin{equation}
	\label{eq:3.12}
		\|D[v]\|_{X_{T,M}}\le\frac{1}{2}\|v\|_{X_{T,M}}
	\end{equation}
	for $v\in X_{T,M}$ and $M\ge M_*$.
	Furthermore, $D[v]$ is bounded and smooth in $\overline{{\mathbb R}^N_+}\times(\tau,T)$ for any $0<\tau<T$.
\end{lemma}
\noindent
{\bf Proof.}
We first prove \eqref{eq:3.12}. 
Let $T>0$, $M\ge 1$, and $1\le q<\infty$.
For $p\in (Nq/(N-1),\infty)$,
by ($G_1$) and \eqref{eq:3.3} we have
\begin{equation}
\label{eq:3.13}
	\begin{aligned}
	\|D[v](t)\|_{L^p} 
 	& 
	\le \int_0^t\|S_1(t-s)F[v](s)\|_{L^p}\,ds\\
 	& 
	\le \int_0^t\|F[v](s)\|_{L^p}\,ds\\
 	& 
	\le C\|v\|_{X_{T,M}}\int_0^te^{Ms}(1+s^{\frac{1}{p}})s^{-\frac{1}{2}}\,ds\\
 	& 
	\le C\|v\|_{X_{T,M}}\,\bigg(1+T^{\frac{1}{p}}\bigg)
	\int_0^te^{Ms}s^{-\frac{1}{2}}\,ds
	\end{aligned}
\end{equation}
for $v\in X_{T,M}$ and $0<t<T$.
Furthermore, 
since it follows from \cite{FIK08}*{Lemma~2.3} that
\begin{equation}
\label{eq:3.14}
	\sup_{x>0}\int_0^\infty\bigg(\frac{|x\pm y|}{t}\bigg)^k\Gamma_1(x\pm y,t)y^{-\frac{j}{2}}\,dy
	\le Ct^{-\frac{k}{2}-\frac{j}{4}},
	\quad k=0,1,\quad j=0,1, 
\end{equation}
for $t>0$,
by \eqref{eq:2.3} and \eqref{eq:3.4} we have
\begin{equation}
\label{eq:3.15}
	\begin{aligned}
 	& 
	|D[v](x,t)|
	\le \int_0^t\int_{\mathbb R^N_+}\Gamma_D(x,y,t-s)|F[v](y,s)|\,dy\,ds
	\\
	& 
	\le C\int_0^t\int_0^\infty\Gamma_1(x_N-y_N,t-s)
	\|F[v](\cdot,y_N,s)\|_{L^\infty({\mathbb R}^{N-1})}\,dy_N\,ds
	\\
	& 
	\le C\|v\|_{X_{T,M}}\int_0^t \int_0^\infty
	\Gamma_1(x_N-y_N,t-s)
	\left(1+(y_N^{-1}s)^{\frac{1}{2}}\right)s^{-\frac{1}{2}}e^{Ms}\,dy_N\,ds
	\\
	& 
	\le C\|v\|_{X_{T,M}}\int_0^t e^{Ms}\left(s^{-\frac{1}{2}}+(t-s)^{-\frac{1}{4}}\right)\,ds
	\end{aligned}
\end{equation}
for $v\in X_{T,M}$, $x\in\overline{{\mathbb R}^N_+}$, and $0<t<T$.
Then,
taking a sufficiently large $M\ge1$ if necessary, 
we can apply Lemma~\ref{Lemma:2.2} to \eqref{eq:3.13} and \eqref{eq:3.15}, 
and for $p\in(Nq/(N-1),\infty]$,
it holds that
\begin{equation}
\label{eq:3.16}
	\sup_{0<t<T}\,e^{-Mt}t^{\frac{N}{2}(\frac{1}{q}-\frac{1}{p})}\|D[v](t)\|_{L^p}
	\le\frac{1}{4}\|v\|_{X_{T,M}}.
\end{equation}

On the other hand, 
we observe from \eqref{eq:2.4} and \eqref{eq:3.4} that
\begin{equation}
\label{eq:3.17}
	\begin{aligned}
	&
	|\partial_{x_N} D[v](x,t)|
	\\
	& \le\int_0^t\int_{{\mathbb R}^N_+}|K(x,y,t-s)||F[v](y,s)|\,dy\,ds\\
	& \le C\int_0^t\int_0^\infty \tilde{K}(x_N,y_N,t-s)
	\|F[v](\cdot,y_N,s)\|_{L^\infty({\mathbb R}^{N-1})}\, dy_N\,ds\\
	&
	\le C\|v\|_{X_{T,M}}
	\int_0^t\int_0^\infty
	e^{Ms}\tilde{K}(x_N,y_N,t-s)
	\bigg(1+(y_N^{-1}s)^{\frac{1}{2}}\bigg)s^{-\frac{1}{2}}\, dy_N\,ds
	\end{aligned}
\end{equation}
for $x\in\overline{{\mathbb R}^N_+}$ and $0<t<T$,
where 
$$
	\tilde{K}(x_N,y_N,t)=\frac{|x_N-y_N|}{t}\Gamma_1(x_N-y_N,t)
	+\frac{x_N+y_N}{t}\Gamma_1(x_N+y_N,t)
$$
for $x_N\ge0$, $y_N>0$ and $t>0$. 
Then, by \eqref{eq:3.14} and \eqref{eq:3.17} we see that
\begin{equation}
\label{eq:3.18}
	\left|\partial_{x_N}D[v](x,t)\right|
	\le C\|v\|_{X_{T,M}}
	\int_0^te^{Ms}\bigg(s^{-\frac{1}{2}}(t-s)^{-\frac{1}{2}}+(t-s)^{-\frac{3}{4}}\bigg)\,ds
\end{equation}
for $x\in\overline{{\mathbb R}^N_+}$ and $0<t<T$.
Furthermore, similarly to \eqref{eq:3.17} and \eqref{eq:3.18}, we see that
$$
	\begin{aligned}
	|\partial_{x_N}D[v](t)|_{L^q}
	&
	\le \int_0^t\int_{\mathbb R^{N-1}}\int_{\mathbb R^N_+}|K(x,y,t-s)||F[v](y,s)|\,dy\, dx'\,ds
	\\
	&
	\le C\int_0^t\int_{\mathbb R^N_+}\tilde{K}(x_N,y_N,t-s)|F[v](y,s)|\,dy\,ds
	\\
	&
	\le C\int_0^t\int_0^\infty\tilde{K}(x_N,y_N,t-s)\|F[v](\cdot,y_N,s)\|_{L^q}\,dy_N\,ds
	\\
	&
	\le C\|v\|_{X_{T,M}}\int_0^t\int_0^\infty e^{Ms}\tilde{K}(x_N,y_N,t-s)
	\bigg(1+(y_N^{-1}s)^{\frac{1}{2}}\bigg)s^{-\frac{1}{2}}\, dy_N\,ds
	\\
	&
	\le C\|v\|_{X_{T,M}}
	\int_0^te^{Ms}\bigg(s^{-\frac{1}{2}}(t-s)^{-\frac{1}{2}}+(t-s)^{-\frac{3}{4}}\bigg)\,ds,
	\qquad\qquad
	t>0.
	\end{aligned}
$$
This together with \eqref{eq:3.18} yields
\begin{equation}
\label{eq:3.19}
	\begin{aligned}
	&
	e^{-Mt}t^{\frac{1}{2}}|\partial_{x_N}D[v](t)|_{L^r}
	\\
	&
	\le C\|v\|_{X_{T,M}}\,e^{-Mt}t^{\frac{1}{2}}
	\int_0^te^{Ms}\bigg(s^{-\frac{1}{2}}(t-s)^{-\frac{1}{2}}+(t-s)^{-\frac{3}{4}}\bigg)\,ds,\qquad r\in[q,p],
	\end{aligned}
\end{equation}
for $0<t<T$.
Moreover, for $p\in(Nq/(N-1),\infty)$,
by \eqref{eq:2.1} and \eqref{eq:3.3} we have
\begin{equation}
\label{eq:3.20}
	\begin{aligned} 
	\|\partial_{x_N}D[v](t)\|_{L^p}
	& \le \int_0^t\|\partial_{x_N}S_1(t-s)F[v](s)\|_{L^p}\,ds
	\\
	& \le C\int_0^t(t-s)^{-\frac{1}{2}}\|F[v](s)\|_{L^p}\,ds
	\\
	&
	\le C\|v\|_{X_{T,M}}\int_0^t(t-s)^{-\frac{1}{2}}(1+s^{\frac{1}{p}})s^{-\frac{1}{2}}e^{Ms}\,ds
	\\
	&
	\le C\|v\|_{X_{T,M}}\,\bigg(1+T^{\frac{1}{p}}\bigg)\int_0^t(t-s)^{-\frac{1}{2}}s^{-\frac{1}{2}}e^{Ms}\,ds
	\end{aligned}
\end{equation}
for $0<t<T$.
Then, by Lemma~\ref{Lemma:2.2} with \eqref{eq:3.18}, \eqref{eq:3.19}, and \eqref{eq:3.20}, 
taking a sufficiently large $M\ge1$ if necessary, 
for $p\in(Nq/(N-1),\infty]$ and $r\in[q,p]$,
we see that 
$$
	\begin{aligned}
	&
	\sup_{0<t<T}\,e^{-Mt}
	t^{\frac{N}{2}(\frac{1}{q}-\frac{1}{p})+\frac{1}{2}}\|\partial_{x_N}D[v](t)\|_{L^p}
	\le\frac{1}{8}\|v\|_{X_{T,M}},
	\\
	&
	\sup_{0<t<T}\,e^{-Mt}t^{\frac{1}{2}}|D[v](t)|_{L^r}
	\le\frac{1}{8}\|v\|_{X_{T,M}}.
	\end{aligned} 
$$
This together with \eqref{eq:3.16} implies \eqref{eq:3.12}.

Next we prove the boundedness and smoothness of $D[v]$.
It follows from the semigroup property of $S_1(t)$ that 
$$
	\begin{aligned}
 	& 
	D[v](x,t)
 	=\int_0^t[S_1(t-s)F[v](s)](x)\,ds
	\\
 	& 
	=S_1(t-\tau)D[v](x,\tau)+\int_{\tau}^t[S_1(t-s)F[v](s)](x)\,ds
	\end{aligned}
$$
for $x\in\overline{{\mathbb R}^N_+}$ and $0<\tau<t<T$.
Then, by \eqref{eq:3.2} and (${\rm G_2}$) 
we see that 
$$
	S_1((t-\tau))D[v](x,\tau)
$$
is bounded and smooth in $\overline{{\mathbb R}^N_+}\times(\tau,T)$. 
Furthermore, by \eqref{eq:3.4} 
we apply the same argument as in \cite{F}*{Section~3, Chapter 1} 
to see that 
$$
	\int_{\tau}^t[S_1(t-s)F[v](s)](x)\,ds
$$
is bounded and smooth in $\overline{{\mathbb R}^N_+}\times(\tau,T)$. 
(See also \cite{FIK01}*{Proposition~5.2} and \cite{IKK01}*{Lemma~2.1}.) 
Then we deduce that 
$D[v]$ and $\partial_{x_N}D[v]$ are bounded and smooth 
in $\overline{{\mathbb R}^N_+}\times(\tau,T)$ for $0<\tau<T$. 
Thus Lemma~\ref{Lemma:3.2} follows.
$\Box$ \vspace{5pt}

Now we are ready to complete the proof of Theorem~\ref{Theorem:1.1}.
\vspace{5pt}
\newline
{\bf Proof of Theorem~\ref{Theorem:1.1}.}
Let $T>0$, $M\ge1$, $1\le q< \infty$, and $p\in(Nq/(N-1),\infty]$
Then, since $\|\varphi\|_{L^r_\alpha}\le \|\varphi\|_{L^r_\beta}$ for  $r\in[1,\infty]$ and $0\le\alpha\le\beta$,
by ($G_1$) and \eqref{eq:2.2} we have
\begin{equation}
\label{eq:3.21}
	e^{-Mt}E[S_1(t)\varphi](t)
	\le (c_1+c_2+c_3)\|\varphi\|_{L^q_{\alpha(p)}}
\end{equation}
for $t>0$,
where $c_1$, $c_2$, and $c_3$ are positive constants given in ($G_1$) and Lemma~\ref{Lemma:2.1}, respectively,
and $\alpha(p)$ is given in \eqref{eq:1.9}.
Furthermore, by Lemma~\ref{Lemma:3.2}, taking a sufficiently large $M\ge1$ if necessary, we see that
\begin{equation}
\label{eq:3.22}
	\|D[v]\|_{X_{T,M}}\le\frac{1}{2}\|v\|_{X_{T,M}},
	\qquad
	v\in X_{T,M},
\end{equation}
for $0<t<T$.
Set
\begin{equation}
\label{eq:3.23}
	m:=2(c_1+c_2+c_3)\|\varphi\|_{L^q_{\alpha(p)}}.
\end{equation}
We deduce from \eqref{eq:3.1}, \eqref{eq:3.21}, \eqref{eq:3.22}, and \eqref{eq:3.23} that
\begin{equation}
\label{eq:3.24}
	\begin{aligned}
	\|Q[v]\|_{X_{T,M}}
	&
	\le
	\sup_{0<t<T}e^{-Mt}E[S_1(t)\varphi](t)+\|D[v]\|_{X_{T,M}}
	\\
	&
	\le (c_1+c_2+c_3)\|\varphi\|_{L^q_{\alpha(p)}}+\frac{1}{2}\|v\|_{X_{T,M}}\le  m
	\end{aligned}
\end{equation}
for $v\in X_{T,M}$ with $\|v\|_{X_{T,M}}\le m$.
Similarly, it follows from \eqref{eq:3.22} that
\begin{equation}
\label{eq:3.25}
	\left\|Q[v_1]-Q[v_2]\right\|_{X_{T,M}}
	=\|D[v_1-v_2]\|_{X_{T,M}}\le\frac{1}{2}\|v_1-v_2\|_{X_{T,M}}
\end{equation}
for $v_i\in X_{T,M}$ $(i=1,2)$. 
Then, by \eqref{eq:3.24} and \eqref{eq:3.25}, 
applying the contraction mapping theorem,
we find a unique solution $v\in X_{T,M}$ with $\|v\|_{X_{T,M}}\le m$ such that
$$
	v=Q[v]=S_1(t)\varphi-D[v](t)
	\quad\mbox{in}\quad X_{T,M}.
$$
In particular, we see that
$$
	\|v\|_{X_{T,M}}\le C\|\varphi\|_{L^q_{\alpha(p)}}.
$$
Moreover,
by ($G_2$) and Lemma~\ref{Lemma:3.2}, 
we see that $v$ is bounded and smooth in $\overline{{\mathbb R}^N_+}\times(T_1,T)$
for any $0<T_1<T$.

Set
$$
w(x,t)=\int_0^t[S_2(t-s)\partial_{x_N}v(s)](x)\, ds
$$
for $x\in\overline{{\mathbb R}^N_+}$ and $t\in(0,T)$. 
By \eqref{eq:2.13} and \eqref{eq:3.23} we obtain
$$
	\begin{aligned}
	\|w(t)\|_{L^p}
	&
	\le \int_0^t\|S_2(t-s)\partial_{x_N}v(s)\|_{L^p}\,ds
	\\
	&
  	\le C\int_0^t\bigg(|\partial_{x_N}v(s)|_{L^q}+|\partial_{x_N}v(s)|_{L^p}\bigg)\,ds
  	\\
  	&
  	\le C\int_0^te^{Ms}s^{-\frac{1}{2}}\|v\|_{X_{T,M}}\,ds
 	\le Ce^{MT}T^{\frac{1}{2}}\|\varphi\|_{L^q_{\alpha(p)}}<\infty,
\end{aligned}
$$
and
$$
	\begin{aligned}
	|w(t)|_{L^r}
	&
	\le \int_0^t|S_2(t-s)\partial_{x_N}v(s)|_{L^r}\,ds
	\\
	&
	\le C\int_0^t|\partial_{x_N}v(s)|_{L^r}\,ds
	\\
	&
 	\le C\int_0^te^{Ms}s^{-\frac{1}{2}}\|v\|_{X_{T,M}}\,ds
 	\le Ce^{MT}T^{\frac{1}{2}}\|\varphi\|_{L^q_{\alpha(p)}}<\infty,
	\end{aligned}
$$
for $0<t<T$.
Furthermore, by ($P_3$) 
we apply an argument similar to that in the proof of Lemma~\ref{Lemma:3.2} 
and see that $w$ is bounded and smooth in $\overline{{\mathbb R}^N_+}\times(T_1,T)$ for any $0<T_1<T$. 
Therefore we deduce that 
$(v,w)$ is a solution of \eqref{eq:1.6} in ${\mathbb R}^N_+\times(0,T)$ 
satisfying \eqref{eq:1.10}. 

Let $(\tilde{v},\tilde{w})$ be a solution of \eqref{eq:1.6} in ${\mathbb R}^N_+\times(0,T_*)$ for any $T_*>T$ and such that $\tilde v\in X_{T_*,M_*}$ with some $M_*>0$.
Then $\tilde v\in X_{T,M}$ and 
since 
$$
	v-\tilde{v}
	=Q[v]-Q[\tilde{v}]
	=D[v-\tilde{v}]\quad\mbox{in}\quad X_{T,M},
$$
by \eqref{eq:3.12} we have 
$$
	\|v-\tilde{v}\|_{X_{T,M}}\le\frac{1}{2}\|v-\tilde{v}\|_{X_{T,M}}.
$$
This implies that $v=\tilde{v}$ in $X_{T,M}$. 
Therefore we deduce that $(v,w)$ is a unique global-in-time solution of \eqref{eq:1.6} 
satisfying \eqref{eq:1.10}.
Thus Theorem~\ref{Theorem:1.1} holds for the case $q<\infty$. 
Furthermore, by \cite{FIK08}*{Theorem~1.1} with $\varepsilon=1$ and $\varphi_b\equiv0$
we have Theorem~\ref{Theorem:1.1} for the case $p=q=\infty$,
and the proof of Theorem~\ref{Theorem:1.1} is complete. 
$\Box$
\vspace{8pt}

\noindent
{\bf Acknowledgment.}
The first author was supported in part by the Slovak
Research and Development Agency under the contract No. APVV-18-0308 and by the VEGA grant
1/0339/21. 
The second and third authors were supported in part by JSPS KAKENHI Grant Number JP19H05599. 
The third author was also supported in part by JSPS KAKENHI Grant Number JP20K03689. 
%

\end{document}